\documentclass[11pt,twoside]{amsart}

\usepackage{pstricks}
\usepackage{multido}
\usepackage{pstcol,pst-text}

\usepackage{graphicx}

\usepackage{amssymb}
\usepackage{amsmath}
\usepackage{amstext}
\usepackage{amsbsy}
\usepackage{amsxtra}
\usepackage{amsfonts}
\usepackage{latexsym}
\usepackage{alltt}
\usepackage{bbm}
 
\theoremstyle{plain}
\newtheorem{theorem}{Theorem}[section]
\newtheorem{corollary}[theorem]{Corollary}
\newtheorem{lemma}[theorem]{Lemma}
\newtheorem{proposition}[theorem]{Proposition}
\newtheorem{remark}[theorem]{Remark}
\newtheorem{definition}[theorem]{Definition}

\newcommand{\C}{\mathbb C}
\newcommand{\R}{\mathbb R}

\newcommand{\Z}{\mathbb Z}

\newcommand{\one}{\mathbbm{1}}
\newcommand{\Ad}{\mathrm{Ad}}

\newcommand{\SL}{\mathrm{SL}(2,\C)}

\newcommand{\SU}{\mathrm{SU}(2)}
\newcommand{\su}{\mathfrak{su}(2)}

\newcommand*{\g}{\mathfrak{g}}
\newcommand*{\gc}{\g^\C}
\renewcommand*{\k}{\mathfrak{k}}
\newcommand*{\m}{\mathfrak{m}}
\newcommand*{\loopE}{\Lambda G_\tau}
\newcommand*{\loopC}{\Lambda G^C_\tau}
\newcommand*{\loopI}{\Lambda_+ G_\tau}
\newcommand*{\potl}{\Lambda_{-1,\infty}\gc}
\newcommand*{\loopc}{\Lambda\gc_\tau}
\newcommand*{\loope}{\Lambda\g_\tau}
\newcommand*{\degone}{\Lambda_1\g_\tau}

\begin{document}
\title{Equivariant harmonic cylinders}

\date{\today}

\author{F.E. Burstall}
\author{M. Kilian}

\address{Mathematical Sciences, 
University of Bath, Bath, BA2 7AY, UK.} 
\email{feb@maths.bath.ac.uk {\rm(Burstall)}}
\email{masmk@maths.bath.ac.uk {\rm(Kilian)}}  

\numberwithin{equation}{section}

\thanks{Kilian supported by EPSRC grant 
GR/S28655/01}

\begin{abstract}
\footnotesize
We prove that a primitive harmonic map 
is equivariant if and only if it admits a holomorphic 
potential of degree one. We investigate when the 
equivariant harmonic map is periodic, and as an 
application discuss constant mean curvature cylinders 
with screw motion symmetries.
\end{abstract}

\maketitle


\begin{center} \sc{Introduction} \end{center}

The theory of harmonic maps, especially those from a Riemann surface
to a Riemannian symmetric space, has been greatly enriched in recent
years by the realisation that they constitute an integrable system
\cite{BurFPP,DorPW,hitchin-harmonic,PinS,TerU,Uhl}.  Thus these maps
admit a spectral deformation (the \emph{associated family});
algebro-geometric (finite type) solutions and an interpretation in
terms of certain holomorphic maps to a loop group.  This last
culminates in the (somewhat indirect) description by
Dorfmeister-Pedit-Wu of all primitive harmonic maps of a simply
connected Riemann surface into a $k$-symmetric space in terms of
\emph{holomorphic potentials}: certain holomorphic $1$-forms with
values in a loop algebra.

In this paper, we study the simplest case of this theory: that of
equivariant harmonic maps where the underlying PDE reduces to an ODE.
We show that such maps are characterised as those with a holomorphic
potential of the simplest kind: the $1$-form has simple poles at zero
and infinity and satisfies a natural reality condition. Along the
way, we show that equivariance is a property which is preserved
under spectral deformation.

A pleasant application of the foregoing theory lies in the fact 
that several types of surface of classical geometric interest are
characterised by harmonicity of an appropriate Gauss map
\cite{Bry:wil,BurHer02,RuhV}.  In particular, the classical theory of
constant mean curvature surfaces in $\R^3$ amounts to the study of
harmonic maps to a $2$-sphere with the link between the loop group
approach and the classical surfaces being provided by the
Sym--Bobenko formula.

We apply our general theory to this case which means that we study
constant mean curvature surfaces with screw motion symmetry.  We find
a very simple proof of a result of Do Carmo--Dajczer \cite{DoCD}
which asserts that these surfaces are precisely those in the
associated family of a Delaunay surface (that is, a constant mean
curvature surface of revolution).  For this, we provide an
interpretation of the Sym--Bobenko formula in terms of a homomorphism
from a loop group to the Euclidean group which may be of independent
interest. 

We then turn to a detailed study of the period problem for constant
mean curvature surfaces with screw motion symmetry.  Armed with our
knowledge of the holomorphic potential, we are able to explicitly
compute, in terms of elliptic functions, the corresponding map into
the loop group (usually, this involves solving a Riemann--Hilbert
problem).  With this in hand, we can prove the existence of
infinitely many non-congruent cylinders in the associated family of
each Delaunay surface.  Otherwise said, in each associated 
family of equivariant harmonic maps $\C \to S^2$ there 
are infinitely many equivariant harmonic tori. By our 
results in the first part, all equivariant harmonic tori 
in the two-sphere arise this way.


\section{Primitive harmonic maps and loop groups}


\subsection{}
We study primitive harmonic maps of a Riemann surface into a
$k$-symmetric space and so begin by recalling the ingredients 
of that story \cite{BurP:dre,BurP_adl}.

Let $G$ be a compact semisimple Lie group. 
A (regular) $k$-symmetric
$G$-space \cite{Kow80} is a coset space $N=G/K$ where
$(G^{\tau})_{0}\subset K\subset G^{\tau}$ for some automorphism
$\tau:G\to G$ of finite order $k\geq 2$.

In particular, a Riemannian symmetric space of compact type is the
same as a $2$-symmetric space.

Let $\g$ be the Lie algebra of $G$.  Then $\tau$ induces an
automorphism, also called $\tau$, of $\g$ whose eigenspace
decomposition gives a $K$-stable $\Z_k$ grading on $\gc$:
$\gc=\sum_{\ell\in\Z_k}\g_\ell$.  Here $\g_\ell$ is the
$\omega^\ell$-eigenspace of $\tau$ for $\omega=e^{2\pi i/k}$.

Define $\m\subset\g$ by $\m^C=\sum_{\ell\neq0}\g_\ell$ to get a
reductive decomposition
\begin{equation}
\label{eq:1}
\g=\k\oplus\m.
\end{equation}

\subsection{} Let $M$ be a Riemann surface. We study maps $\varphi:M\to
G/K$ via their \emph{frames}: that is, maps $F:M\to G$ for which
$\varphi=FK$.  Given a frame $F$, set $\alpha=F^{-1}dF$, a $\g$-valued
$1$-form on $M$ and write $\alpha=\alpha_\k+\alpha_\m$ according to
the decomposition \eqref{eq:1}.  Further write
$\alpha_\m=\alpha'_\m+\alpha''_\m$ according to the type
decomposition $TM^\C=T^{1,0}\oplus T^{0,1}$.  Thus $\alpha'_\m$ is a
$(1,0)$-form with values in $\m^\C$ and
$\alpha''_\m=\overline{\alpha'_\m}$.

A map $\varphi:M\to G/K$ is \emph{harmonic} if, for any (hence every) frame,
\begin{equation*}
d\alpha'_\m+[\alpha_k\wedge\alpha'_\m]=0
\end{equation*}
and \emph{primitive} if $\alpha'_m$ takes values in $\g_{-1}$.  It is
easy to see that when $k>2$, a primitive map is automatically
harmonic while, when $k=2$, the primitivity condition is vacuous.  We
combine the cases of interest to us by saying that $\varphi$ is
\emph{primitive harmonic} if $k=2$ and $\varphi$ is harmonic or $k>2$
and $\varphi$ is primitive.

\subsection{} The basic observation \cite{Poh,Uhl,ZakS2} is that if $F$ frames
a primitive harmonic map then defining
\begin{equation}
\label{eq:3}
\alpha_\lambda=\alpha_\k+\lambda^{-1}\alpha'_\m+\lambda\alpha''_\m
\end{equation}
yields a solution of the Maurer--Cartan equations for each non--zero 
complex number $\lambda\in\C^{\times}$. 
Thus, we can (locally or after passage to
the universal cover of $M$) integrate to find
$F_\lambda:M\to G^\C$ such that $F_\lambda^{-1}
dF_\lambda=\alpha_\lambda$.  It is easy to see that $\alpha_\lambda$
has the following symmetries:
\begin{equation*}
\overline{\alpha_{\bar{\lambda}}}=\alpha_\lambda\qquad
\tau\alpha_\lambda=\alpha_{\omega\lambda}.
\end{equation*}
Thus we may choose the constants of integration to ensure:
\begin{enumerate}
\item $\overline{F_{1/\bar\lambda}}=F_\lambda$, for all
$\lambda\in\C^\times$.  Here the conjugation on $G^\C$ has fixed set
$G$ so that, in particular, $F_\lambda$ takes
values in $G$ when $\lambda\in S^1$;
\item $\tau F_\lambda=F_{\omega\lambda}$, for all
$\lambda\in\C^\times$;
\item $F_1=F$;
\end{enumerate}
We encapsulate all this by viewing $F_\lambda$ as a map from (the
universal cover of) $M$ into
the loop group $\loopE$ given by
\begin{equation*}
\loopE=\left\{\text{$g:S^1\to G$ smooth}: g(\omega\lambda)=\tau
  g(\lambda)\right\} 
\end{equation*}
and say that $F_\lambda:M\to\loopE$ is an \emph{extended frame} if
$\alpha_\lambda=F_\lambda^{-1}dF_{\lambda}$ is of the form
\eqref{eq:3}.

We have just seen that any primitive harmonic map has an extended
frame with $F_1=F$.  Conversely, given an extended frame $F_\lambda$,
it is easy to see that, for each $\lambda\in S^1$, $F_\lambda$ frames
a primitive harmonic map: these constitute the \emph{associated
family} of such maps.

\subsection{}
The loop group $\loopE$ participates in an Iwasawa decomposition of
its complexification:  set
\begin{align*}
\loopC&=\left\{\text{$g:S^1\to G^\C$ smooth}: g(\omega\lambda)=\tau
  g(\lambda)\right\}\\
\loopI&=\left\{ \text{$g \in \loopC$: $g$ extends holomorphically to
  $|\lambda|<1$ and $g(0)\in B$}\right\}
\end{align*}
where $K^\C=KB$ is some fixed Iwawsawa decomposition of $K^\C$.
\begin{theorem}[\cite{DorPW,PreS}]\label{h:Iwa}
Pointwise multiplication $\loopE \times \loopI\to\loopC$ is a
diffeomorphism.
\end{theorem}
Consequently, every $g\in\loopC$ can be uniquely factored $g = F\,b$
with $F\in\loopE$ and $b \in \loopI$.

\subsection{}
Theorem~\ref{h:Iwa} underlies a construction of extended frames which
is an infinite-dimensional version of a formula of Symes \cite{Sym80} from
integrable systems theory.  For this, let $\loope,\loopc$ be the Lie algebras
of $\loopE,\loopC$ respectively.  Thus
\begin{align*}
\loopc&=\left\{\text{$\xi:S^1\to \gc$ smooth}: \xi(\omega\lambda)=\tau
  \xi(\lambda)\right\};\\
\loopc&=\left\{\xi\in\loopc: \xi:S^1\to\g\right\}.
\end{align*}
Further let $\potl$ be the vector subspace of $\loopc$ given by
\begin{equation*}
\potl=\left\{\xi\in\loopc:\text{$\xi$ extends holomorphically to
  $0<\vert\lambda\vert<1$ with a simple pole at $0$}\right\}.
\end{equation*}
We have:
\begin{theorem}[\cite{BurP:dre}]\label{h:Symes}
For $z\in\C$ and $\xi\in\potl$, let $\exp(z\,\xi)=F_\lambda(z)b(z)$ be
the Iwasawa decomposition of $\exp(z\,\xi)$.  Then
$F_\lambda:\C=\R^2\to\loopE$ is an extended frame.
\end{theorem}

This construction accounts for all primitive harmonic maps of
semisimple finite type \cite{BurP:dre} and is a special case of the
method of Dorfmeister--Pedit--Wu \cite{DorPW} which describes all
primitive harmonic maps of a simply connected Riemann surface into a
$k$-symmetric space.  Following \cite{DorPW}, we say that the
\emph{holomorphic potential} $\xi dz$ \emph{generates} the associated
family of primitive harmonic maps.


\section{Equivariant primitive harmonic maps}


\subsection{}\label{sec:equi}
Consider a map 
$\varphi:\R^2 \to G/K$ that is 
$\R$-equivariant, so that 
\begin{equation} \label{eq:equivariance}
  \varphi(x,\,y) = 
  \exp(x\,A_0)\,\psi(y)
\end{equation}
for some $A_0 \in \mathfrak{g}$ and map $\psi:\R \to G/K$.
If $g$ frames $\psi$, then $F : \R^2 \to G$ defined by 
$F(x,\,y) = \exp(x\,A_0)\,g(y)$ frames $\varphi$ and 
\begin{equation} \begin{split}
  F^{-1}dF &= \Ad g^{-1}A_0\,dx + g^{-1}g'\,dy \\
           &= A(y)\,dx + B(y)\,dy. \label{eq:y-form}
           \end{split}
\end{equation}
In particular, both components of $F^{-1}dF$ depend 
only on $y$. 

Conversely, given $F:\R^2 \to G$ with $F^{-1}dF$ of the 
form \eqref{eq:y-form}, set 
$g(y) := F(0,\,y)$ and $h(x,\,y) := F(x,\,y)\,g^{-1}(y)$, 
and note the following properties: 
\begin{enumerate}
  \item $g^{-1}dg = B\,dy$
  \item $h(0,\,y) = \one$ for all $y \in \R$
  \item $h^{-1}dh = \Ad g~(F^{-1}dF -g^{-1}dg) = \Ad g~A\,dx$
\end{enumerate}
The Maurer--Cartan equation for $h^{-1}dh$ reads
$\partial/\partial y (\Ad g~A) = 0$. Thus 
$\Ad g~A =: A_0$ is constant, which with (ii) above 
yields $h(x,\,y) = \exp(x\,A_0)$. 
Hence $F(x,\,y) = \exp(x\,A_0)\,g(y)$, so that $F$ frames 
an equivariant map. This proves
\begin{proposition}\label{th:equi_y}
$F(x,\,y):\R^2 \to G$ frames an equivariant map if and only if 
$F^{-1}dF$ depends only on $y$.
\end{proposition}


\subsection{}\label{ss:degree1}
We now come to the first of our main results: 
we consider equivariant
primitive harmonic maps and show that these are all generated by
holomorphic potentials of the simplest kind.  
Moreover, we find that
equivariance is shared by all members of an associated family.

For this, contemplate the \emph{degree one} elements of $\loope$:
\begin{equation*}
\degone=\left\{\xi\in\loope:
  \xi=\lambda\xi_1+\xi_0+\lambda^{-1}\xi_{-1}\right\}=\potl\cap\loope.
\end{equation*}
Note that $\xi_\ell\in\g_\ell$ and $\overline{\xi_\ell}=\xi_{-\ell}$. 
We call a holomorphic potential $\xi dz$ with $\xi\in\degone$ a 
\emph{degree one potential}.

Given a degree one potential $\xi dz$ with $\xi\in\degone$, we 
examine the extended frame $F_\lambda:\R^2\to\loopE$ that results
from Theorem~\ref{h:Symes}.  With $z=x+iy$, we have
$\exp(z\xi)=\exp(x\xi)\exp(iy\xi)$ and since $\exp(x\xi)\in\loopE$,
for all $x\in\R$, we conclude that
\begin{equation} \label{eq:F_split}
  F_\lambda (z) = \exp(x\,\xi)\,g_\lambda (y)
\end{equation}
for some $g_\lambda:\R \to \loopE$ resulting from the Iwasawa
decomposition of $y\mapsto\exp(iy\xi)$.  Hence $F_\lambda$ is the
frame of an equivariant primitive harmonic map for each $\lambda \in
S^1$.  We have therefore proved
\begin{proposition} \label{th:collin1}
A degree one potential generates an associated 
family of equivariant primitive harmonic maps.
\end{proposition} 

\subsection{} We next prove the converse of Proposition
\ref{th:collin1}, that any equivariant primitive 
harmonic map is generated by a
degree one potential.  So let $\varphi:\R^2 \to G/K$ be an
equivariant primitive harmonic map.  After translation by an element
of $G$, we may assume that $\varphi(0)=eK$ and then we can find an
equivariant frame $F$ as in paragraph~\ref{sec:equi} with $F(0) = \one$.

From Proposition \ref{th:equi_y} we know that $\alpha = F^{-1}dF$ is
of the form $\alpha = A(y)\,dx + B(y)\,dy$, so that $\alpha_\lambda =
A_\lambda(y)\,dx + B_\lambda(y)\,dy$ with
\begin{itemize}
\item[(a)] $A_\lambda(y),\,B_\lambda(y) \in 
  \Lambda_1\,\mathfrak{g}_\tau$,
\item[(b)] $\alpha_\lambda(\partial/\partial \bar{z})$ 
  is holomorphic (indeed affine) in $\lambda$ on $\C$.
\end{itemize}
Let $F_\lambda$ be the corresponding extended frame.  An immediate
consequence of (a) and Proposition \ref{th:equi_y} is that, for each
$\lambda \in S^1$, $F_\lambda$ frames an equivariant map so that we
may conclude:
\begin{proposition}
If a primitive harmonic map is equivariant, then so are all 
members of its associated family.
\end{proposition}
%

\subsection{} To show that an extended frame 
$F_\lambda$ of an associated family of equivariant 
primitive harmonic maps is generated by a degree one potential, 
we argue as in the proof of Proposition \ref{th:equi_y}: 
Set 
$g_\lambda(y) = F_\lambda(0,\,y)$ 
and conclude that $\Ad g_\lambda~A_\lambda$ is constant. 
Since $g_\lambda(0) = \one$, that constant is 
$A_\lambda(0) \in \Lambda_1\mathfrak{g}_\tau$. Set
\[
  \xi = A_\lambda(0) = \Ad g_\lambda A_\lambda(y),
\]
to obtain $F_\lambda(x,\,y) = \exp(x\,\xi)\,g_\lambda(y)$. 
It remains to show that $F_\lambda$ arises from the 
holomorphic potential $\xi dz$, that is 
$\exp(z\,\xi) = F_\lambda\,b_\lambda$ for some 
smooth map $b_\lambda:\R \to \Lambda_+ G^\C_\tau$. 
For this it suffices to check that 
$b_\lambda(y) = g^{-1}_\lambda(y)\,\exp(iy\,\xi)$ 
is holomorphic in $\lambda$ near $0$. This is certainly 
true at $z=0$, so we need only check holomorphicity in 
$\lambda$ of its right(!) Maurer--Cartan form
\begin{equation*}
  db_\lambda b^{-1}_\lambda =(\Ad g^{-1}~i\xi-g^{-1}g')~dy
  = i(A_\lambda +i B_\lambda)~dy 
  = \tfrac{i}{2}\alpha_\lambda(\partial/\partial\bar{z}),
\end{equation*}
which is clearly holomorphic in $\lambda$ by (b) above. 
Thus $F_\lambda(z)$ is the $\Lambda G_\tau$-factor of 
$\exp(z\,\xi)$, and we have proven that equivariant 
harmonic maps are generated by degree one potentials. 
Together with Proposition \ref{th:collin1}, we have proven
\begin{theorem} \label{th:Collins} 
A primitive harmonic map is equivariant if and only if it 
is generated by a degree one potential.
\end{theorem}
%


\section{The Euclidean group and the Sym--Bobenko formula}


\subsection{} 
Consider the semi-direct product $G \ltimes \mathfrak{g}$ where 
$G$ acts on $\mathfrak{g}$ via the adjoint action. Thus
\[
  (g,\,\zeta)\,(h,\,\eta) = (gh,\,\Ad g~\eta + \zeta).
\]
$G \ltimes \mathfrak{g}$ has an affine action on $\mathfrak{g}$ via
\begin{equation} \label{eq:aff_act}
  (g,\,\zeta) \cdot \eta = \Ad g~\eta + \zeta.
\end{equation}

\begin{proposition}
  For $\mu \in S^1$ we have a homomorphism
\begin{equation}
  \Phi_\mu : \Lambda G_\tau \to G \ltimes \mathfrak{g},\,\,
  ~F \mapsto \left.\left(F,\,F'
    F^{-1} \right) \right|_{\mu}
\end{equation}
where, for $F \in \Lambda G_\tau$ and $\lambda = e^{it}$, $F'$ denotes the 
derivative of $F$ with respect to $t \in \R$.
\end{proposition}
\begin{proof}
Let $F_1,\,F_2 \in \Lambda G_\tau$. Then
\begin{equation*}
\Phi_\mu (F_1F_2) = 
  \left.\left(F_1F_2,\,\Ad F_1 (F_2'F_2^{-1}) +  
      F_1'   F_1^{-1}\right) \right|_{\mu}  
  = \Phi_\mu (F_1) \,\Phi_\mu (F_2).
\end{equation*}
\end{proof}
Differentiating $\Phi_\mu$ at $\one \in \Lambda G_\tau$ provides a 
Lie algebra homomorphism, also denoted by 
$\Phi_\mu: \Lambda \mathfrak{g}_\tau \to 
\mathfrak{g} \ltimes \mathfrak{g}$ and given by
\begin{equation*}
  \Phi_\mu (\xi) = \left.\left(\xi,\,\xi'\right) \right|_{\mu},
\end{equation*}
where we again write $\xi' = (\partial/ \partial t)\,\xi$, for $\xi
\in \Lambda \mathfrak{g}_\tau$. 
%


\subsection{}
All this comes alive when $G=\SU$. Then $\su \cong \R^3$ and the 
affine action \eqref{eq:aff_act} gives a double cover 
$\SU \ltimes \su \to \mathrm{Euc}(3)$ of the Euclidean group.
An extended frame $F(z,\,\lambda)$ gives rise to an associated 
family of parallel pairs $f^{\pm}_\lambda$ of 
conformal constant mean curvature immersions via the Sym--Bobenko
formula \cite{Bob:cmc,Sym} which, in our formalism, reads 
\begin{equation} \label{eq:sym-bob} 
  f^{\pm}_\lambda = -\tfrac{1}{2H}\,\Phi_\lambda(F)\cdot \left( \mp 
  \mathrm{e}_1 \right).
\end{equation}
Here $\mathrm{e}_1 \in \su$ is the normal to $f=f^+_1$ at $z=0$, 
and $H \in \R^*$ is the mean curvature. Since the normal to 
the surface $f^-$ in \eqref{eq:sym-bob} is 
$N_\lambda = \mathrm{Ad}F_\lambda \,\mathrm{e}_1$, 
the parallel surfaces satisfy
\begin{equation}
  f_\lambda^- + \tfrac{1}{H} \,N_\lambda = f_\lambda^+.
\end{equation}
Suppose that $F_\lambda(z)$ is generated by 
$\xi dz$ with $\xi \in \Lambda_1 \su_\tau$. Then, as  
in \eqref{eq:F_split} we have that 
$F_\lambda(x,\,y) = \exp(x\xi)\,g_\lambda(y)$, and 
the associated families are 
\begin{align*}
  f^{\pm}_\lambda(x,\,y) &= -\tfrac{1}{H}\,\Phi_\lambda\left( 
    \exp(x\xi) \right) \, \Phi_\lambda (g_\lambda(y)) 
  \cdot \left( \mp\,\tfrac{1}{2}\,\mathrm{e}_1 \right)\\
  &= -\tfrac{1}{H}\, 
    \exp(x\Phi_\lambda (\xi)) \, \cdot (-H\,f_\lambda^{\pm}(0,\,y)).
\end{align*}
Thus at $\lambda \in S^1$ the surfaces $f^{\pm}_\lambda$ 
have screw-motion symmetry generated by 
$\Phi_\lambda(\xi) \in \mathfrak{euc}(3) = \su \ltimes \su$, and  
we have proven
\begin{proposition} \label{th:screw} A degree one potential $\xi dz$
generates associated families $f^{\pm}_\lambda$ of constant mean
curvature surfaces, of which each member has screw-motion symmetry
generated by $\Phi_\lambda(\xi)$.
\end{proposition}
We will show that at specific $\mu\in S^1$ we obtain 
surfaces of revolution and thus Delaunay surfaces. This 
will be the case exactly when the orbits of the 
$1$-parameter subgroup generated by $\Phi_\mu(\xi)$ 
are co-axial circles.
\begin{lemma} \label{th:perp}
An element $(\xi,\,\eta) \in \su \ltimes \su$ generates 
a rotation about an axis if and only if $\xi \perp \eta$ 
with respect to the Killing form.
\end{lemma}
\begin{proof}
Let $(\xi,\,\eta) \in \su \ltimes \su$ and write 
$\eta = \eta^\top + \eta^\bot$ with 
$\eta^\top \in \mathrm{Im}(\mathrm{ad}\,\xi)$ and 
$\eta^\top \perp \eta^\bot$ 
( and thus $\eta^\bot \parallel \xi$ ). Hence there exists 
a unique $\zeta \in \mathrm{Im}(\mathrm{ad}\,\xi)$ with 
$[\xi,\,\zeta]=-\eta^\top$. 
The adjoint action of $\SU \ltimes \su$ reads 
\begin{equation*}
  \Ad (g,\,\zeta)\,(\xi,\,\eta) = 
  (\Ad g~\xi,\,\Ad g~\eta-[\Ad g~\xi,\,\eta]),
\end{equation*}
and we compute
\begin{align*}
  (\one,\,\zeta)^{-1} \exp\left(t\,(\xi,\,\eta)\right)\,
  (\one,\,\zeta) &= \exp \left( t \,\Ad (\one,\,-\zeta)
    (\xi,\,\eta)\right) 
  = \exp \left( t\,(\xi,\,\eta + [\xi,\,\zeta])\,\right) \\
  &= \exp \left( t\,(\xi,\,\eta^\bot)\,\right) 
  = (\exp(t\,\xi),\,t\,\eta^\bot),
\end{align*}
since $[\xi,\,\eta^\bot]=0$. Thus $(\one,\,\zeta)^{-1} 
\exp\left(t\,(\xi,\,\eta)\right)\,(\one,\,\zeta)$ is a 
screw-motion symmetry with axis $\langle \xi \rangle$ 
and translational part given by $\eta^\bot$. It follows 
that $\exp(t(\xi,\,\eta))$ is a screw-motion symmetry 
with axis $\langle \xi \rangle + \zeta$ and the same 
translational part. In particular, if $\xi \perp \eta$ 
and thus $\eta^\bot=0$, then $\exp(t(\xi,\,\eta))$ consists 
of rotations in $\xi^\bot$ with axis 
$\langle \xi \rangle + \zeta$.
\end{proof}

\begin{theorem} \label{th:DEL}
A degree one potential $\xi dz$ with 
$\xi \in \Lambda_1\su_\tau$ 
generates an associated family of a Delaunay surface.
\end{theorem}
\begin{proof}
Write $\xi = \lambda^{-1}\xi_- + \xi_0 + \lambda\,\xi_+$ 
so that 
$\Phi_\mu(\xi) = (\mu^{-1}\xi_- + \xi_0 + \mu\,\xi_+,\,
i\mu\,\xi_+ - i\mu^{-1}\xi_-)$. Now $\xi_{\pm} \perp \xi_0$ 
since $\g_\pm\perp\g_0$, so with $\langle \,,\, \rangle$ 
denoting the Killing form, we have
$\langle \mu^{-1}\xi_- + \xi_0 + \mu\,\xi_+,\,
  i\mu\,\xi_+ - i\mu^{-1}\xi_- \rangle = 
  i\mu^2\langle \xi_+,\,\xi_+ \rangle - 
  i\mu^{-2}\langle \xi_-,\,\xi_- \rangle$, which vanishes for 
$\mu^4 = \langle \xi_-,\,\xi_- \rangle/
\langle \xi_+,\,\xi_+ \rangle$.
Clearly the fourth roots are on $S^1$ since 
$\overline{\xi_-} = \xi_+$. 
Hence at these points, by Lemma~\ref{th:perp}, the surface is 
a surface of revolution and thus a Delaunay surface.
Hence a degree one potential generates an associated family 
of a Delaunay surface.
\end{proof}
With the above results we are now in a position to present 
a new and elementary proof of a theorem by DoCarmo and 
Dajczer \cite{DoCD}. 
\begin{theorem}\cite{DoCD} 
A constant mean curvature surface has screw-motion symmetry 
if and only if it lies in an associated 
family of a Delaunay surface.
\end{theorem}
\begin{proof}
Let $y \mapsto (u(y),\,v(y),\,w(y))^t$ be the generating 
curve in $\R^3$ of a helicoidal surface such that we have 
a conformal parametrisation 
\begin{equation}
f(x,\,y) = \begin{pmatrix} \cos(x) & \sin(x) & 0 \\
  -\sin(x) & \cos(x) & 0 \\ 
  0 & 0 & 1 \end{pmatrix}
  \begin{pmatrix} u(y) \\ v(y) \\ w(y) \end{pmatrix} 
  + \begin{pmatrix} 0 \\ 0 \\ a\,x \end{pmatrix},\,a \in \R. 
\end{equation}
The Gauss map $N = |f_x \times f_y|^{-1}f_x \times f_y$ is 
equivariant, since 
\begin{equation*}
f_x \times f_y =  \begin{pmatrix} 
        \cos(x) & \sin(x) & 0 \\
        -\sin(x) & \cos(x) & 0 \\
        0 & 0 & 1 \end{pmatrix}\,
        \begin{pmatrix}
        -uw'-av' \\ -vw'+au' \\ uu'+vv' \end{pmatrix},
\end{equation*}
and $|f_x \times f_y|$ is independent of the parameter $x$. 
Since the mean curvature is constant, the Gauss map is 
harmonic \cite{RuhV}, and hence by Theorem \ref{th:Collins} 
admits a holomorphic potential $\xi dz$ with 
$\xi \in \Lambda_1 \su_\tau$ of degree one. 
By Theorem \ref{th:DEL} this generates the harmonic map of an 
associated family of Delaunay surfaces.

Conversely, let $f_{\lambda_0}:\R^2 \to \R^3$ 
be some member of an associated family of 
a Delaunay surface. Then at some $\lambda_1 \in S^1$ we 
have that $f_{\lambda_1}$ is a surface of revolution, 
and thus has an equivariant Gauss map generated 
by some degree one potential. 
By Proposition \ref{th:screw}, the 
surface has screw-motion symmetry.
\end{proof}
%


\section{A normal form}


Our main objective in this second part 
is to solve period problems for equivariant 
harmonic maps $\C \to S^2$ and the resulting 
helicoidal constant mean curvature surfaces. 
We proceed by deriving a normal form 
for the corresponding potentials of degree one, 
and then compute the surface invariants generated by 
such a potential. These are necessary ingredients for 
the subsequent explicit extended frame and its monodromy. 

We saw above that all equivariant 
harmonic maps $\C \to S^2$ come from holomorphic degree one 
potentials on $\C$ with values in $\Lambda_1\su_\tau$. 
Such a potential is of the form 
\begin{equation} \label{eq:unnormal}
  \begin{pmatrix} c & a\,\lambda^{-1} + \bar{b}\,\lambda \\ 
   b\,\lambda^{-1} + \bar{a}\,\lambda & -c \end{pmatrix}\,i\,dz
\end{equation}
with $a,\,b \in \C$ and $c \in \R$. It turns out that 
up to rigid motions and a translation of $y$, all equivariant harmonic maps 
$\C \to S^2$ can be generated by a two-parameter 
subfamily of such potentials. This prompts the following 
\begin{definition}
A $\Lambda_1\su_\tau$-valued holomorphic $1$-form 
on $\C$ of the form
\begin{equation}\label{eq:del_mod}
  \begin{pmatrix} 0 & a\,\lambda^{-1} + b\,\lambda \\ 
   b\,\lambda^{-1} + a\,\lambda & 0 \end{pmatrix}\,i\,dz 
\end{equation}
with $a,\,b \in \R$, is said to be in \emph{normal form}.
\end{definition}
\begin{lemma} \label{th:normalform}
Up to rigid motions, any Delaunay surface 
can be generated by a potential in normal form.
\begin{proof}
Let $f:\R^2 \to \R^3$ be a conformal immersion 
of a Delaunay surface. 
Then its Gauss map $N$ is equivariant and thus 
generated by a degree one potential $\eta\, dz$ with 
$\eta \in \Lambda_1\su_\tau$ such that for the associated 
family of Gauss maps $N_\lambda$ we have $N_1=N$. Then 
up to rigid motion we have for the associated family of 
surfaces $f_\lambda$ that $f_1=f$. 
Apriori, $\eta\,dz$ is of the form \eqref{eq:unnormal}
with $a,\,b \in \C$ and $c \in \R$. 
Since $f_1$ is a surface of revolution, 
Lemma \ref{th:perp} implies that $\eta \perp \eta'$ at 
$\lambda = 1$. 
A computation reveals that 
$\langle \eta,\,\eta'\rangle|_{\lambda=1} =0$ 
if and only if the product $ab \in \R$.

After a diagonal unitary gauge (rotation in the tangent plane), 
we may assume that $a \in \R$, and consequently also $b \in \R$. 
Finally, after a possible translation in the domain, 
we may assume that 
the base point lies on a circle of maximal or minimal radius. 
This condition 
is equivalent to $N_1(0) = \mathbf{e}_1$ being perpendicular 
to the axis of revolution $\eta_{\lambda=1}$. A quick 
computation shows that this holds if and only if $c=0$. 
\end{proof}
\end{lemma}
\begin{corollary} \label{cor:equi_harm}
Up to isometry, any equivariant harmonic map 
$\C \to S^2$ is generated by a $\Lambda_1\su_\tau$-valued 
holomorphic $1$-form in normal form.
\end{corollary}
We will compute the metric and Hopf differential of 
the associated Delaunay surface generated by a potential in 
normal form. 
Since the parameter in all further elliptic functions 
and integrals is in terms of $a^2/b^2$, we set
\begin{equation} \label{eq:kappa}
  \kappa = a^2/b^2 \mbox{ and } \kappa' = 1- \kappa.
\end{equation}
\begin{theorem} \label{th:invariants}
Let $H\in \R^*$ and consider an associated family $f_\lambda$ 
generated by a potential in normal form with constants 
$a,\,b \in \R$ via the Sym--Bobenko formula \eqref{eq:sym-bob}.
The Hopf differential of $f_\lambda$ 
is $\lambda^{-2}\,Q\,dz^2$ with 
\begin{equation} \label{eq:hopf}
 Q = 2\,a\,b\, H^{-1}.
\end{equation}
The metric is given by $v_0^2\,(dz \otimes d\bar{z})$ with 
$v_0 = v_0(y)$ the Jacobian elliptic function 
\begin{equation} \label{eq:v_Delaunay}
   v_0(y) = 2bH^{-1}\,\mathrm{dn}(2by\,|\,\kappa').
\end{equation}
\end{theorem}
\begin{proof} 
Let $F : \R^2 \to \Lambda \SU_\tau$ be the extended framing 
obtained from a potential $\xi dz$ in normal form and let $f$ 
be one of the immersions given by \eqref{eq:sym-bob}. 
This means that $F$ is the 
$\Lambda \SU_\tau$-factor of $\exp(z\,\xi) = F\,B$. 
Then $\alpha = F^{-1}dF$ is of the form
\begin{equation}\label{eq:MC_loop} 
  \alpha = \frac{1}{2v}
  \begin{pmatrix} - v_zdz + v_{\bar{z}} d\bar{z} & 
    2i\lambda^{-1}Q\,dz + iv^2 \lambda H\,d\bar{z} \\ 
    iv^2\lambda^{-1}Hdz + 2i\lambda \overline{Q}\,d\bar{z} & 
    v_zdz - v_{\bar{z}} d\bar{z} 
  \end{pmatrix} 
\end{equation}
for functions $v,\,Q,\,H$. 
The smooth function $v^2:\R^2 \to \R_+$ is the conformal factor, 
the quadratic differential $\lambda^{-2}Q\,dz^2$ the 
Hopf differential of $f_\lambda$.

Since $\xi\,\,dz = B^{-1}\,\alpha \, B + B^{-1}\,d\,B$, 
comparisson of $\lambda^{-1}$ coefficients yields 
$\alpha_{-1} = B_0\,\xi_{-1} \, B_0^{-1}$. 
Writing 
\begin{equation*}
  B_0 = \begin{pmatrix} r & 0 \\ 0 & 1/r \end{pmatrix} 
\end{equation*}
for smooth $r:\R \to \R_+$ 
this gives $i v^{-1} Q = i r^2 a$ and $iv\, H = 2r^{-2} b$. 
Eliminating $v$ from these two equations yields the formula 
for the Hopf differential \eqref{eq:hopf}.
Furthermore, the second of these equations evaluated at the 
basepoint $z=0$ reads $v(0)\,H = 2r^{-2}(0)\,b$. 
{}From $F(0) = B(0) = \one$ it follows that $r(0) = 1$ and 
consequently the function $v:\R \to \R^*$ generated by the 
potential $\xi\,\,dz$ satisfies
\begin{equation}\label{eq:v(0)}
  v(0) = 2\,b/H.
\end{equation}
Since the extended frame splits as in \eqref{eq:F_split}, 
the conformal factor $v^2 = v^2(y)$ depends only on the 
variable $y \in \R$, and the Gau{\ss} equation is 
\begin{equation} \label{eq:gauss_x}
  v^{-1}\ddot{v} - v^{-2} \dot{v}^2 +
  v^2 H^2 - 4\,v^{-2}|Q|^2 = 0.
\end{equation}
Since $H$ and $Q$ are constant, we can find a first integral of 
equation \eqref{eq:gauss_x} of the form 
$\dot{v}^2 = A\,(v^2 - B)\,(v^2 - C)$. Differentiating this, 
and using $A(B+C) = A\,v^2-v^{-2}\dot{v}^2 + v^{-2}ABC$, gives 
$v^{-1}\ddot{v} - v^{-2} \dot{v}^2 - A\, v^2  + v^{-2}ABC = 0$.
Comparing this with \eqref{eq:gauss_x} implies that 
the constants $A,\,B$ and $C$ must satisfy 
$A = -H^2 \mbox{ and } BC = 4\,H^{-2} |Q|^2$. 
Set $B = 4H^{-2}a^2$ and $C = 4H^{-2}b^2$. Hence 
\begin{equation} \label{eq:v_elliptic} 
  \dot{v}^2 = - H^2\bigl( v^2 - 4H^{-2}a^2 \, \bigr)
  \,\bigl( v^2 - 4H^{-2}b^2\,\bigr) 
\end{equation}
\begin{remark} \label{rem:v_elliptic}
Two properties follow from the
differential equation \eqref{eq:v_elliptic}: 

\noindent (i) $\dot{v}(y)=0$ 
if and only if $v(y) \in \{ \pm 2a/H,\,\pm 2b/H \}$. 

\noindent (ii) Since the left hand side in 
\eqref{eq:v_elliptic} is positive, the 
signs of the two factors on the right 
hand side in \eqref{eq:v_elliptic} must be different, 
and thus any solution $v$ oscillates either between 
$-2|a/H|$ and $-2|b/H|$, or between $2|a/H|$ and $2|b/H|$. 
\end{remark}
The general solution of \eqref{eq:gauss_x} is given 
in terms of a Jacobi elliptic function as follows: 
Taking the square root, rewrite \eqref{eq:v_elliptic} as 
\begin{equation} \label{eq:sign}
  dy = \pm\,\left(- H^2(v^2-4H^{-2}a^2)
    (v^2-4H^{-2}b^2)\right)^{-1/2}dv.
\end{equation}
Integrating, we pick up a constant $C \in \R$, 
and substituting $2|a|\,t = v\,|H|$ gives 
\begin{equation*} 
    y = C \, \pm \,\frac{1}{2i|b|}\int_0^{\tfrac{|H|\,v}{2|a|}} 
      \hspace{-5mm} \frac{dt}{\sqrt{(1-t^2)\,
          (1 - \kappa t^2)}} 
      = C \, \pm \,\frac{1}{2i|b|} \,\mathrm{F}\, 
      \left( \left. \arcsin(\tfrac{|H|\,v}{2|a|})\,\right|
      \,\kappa \,\right),
\end{equation*}
where $\mathrm{F}$ denotes the elliptic integral 
of the first kind 
\begin{equation} \label{eq:F}
  \mathrm{F}(\varphi \,|\,m) =  
  \int_0^{\sin \varphi} \hspace{-4mm}
  \frac{dt}{\sqrt{(1-t^2)\,(1-m\,t^2)}} .
\end{equation}
Thus the general solution of \eqref{eq:gauss_x} is 
\begin{equation}
  v(y) = 2|a/H|\,\mathrm{sn}( \pm 2i|b|(y - C) 
  \left.\right|\,\kappa \,).
\end{equation} 
It remains to determine the solution $v_0$ with 
initial condition $v(0) = 2b/H$. 
The solution with $v(0) = 2|b/H|$ is given by 
$v(y) = \pm 2|a/H|\,\mathrm{sn}( 2i|b|y + 
  \left. \mathrm{F} (\arcsin (\kappa^{-1/2})\,|\,\kappa) 
  \,\right|\,\kappa)$. 
By the \emph{complex arguments} formula for $\mathrm{sn}$, 
see 16.21.1 in \cite{AbrS}, this simplifies to 
$v(y) = \pm 2|b/H|\,\mathrm{dn}(2|b|y\,|\,\kappa')$. 
Choosing the sign in \eqref{eq:sign} so that 
$\pm |b/H| = b/H$, 
and using the fact that $\mathrm{dn}$ is an even 
function 16.8.3 \cite{AbrS}, proves \eqref{eq:v_Delaunay} 
and concludes the proof of Theorem \ref{th:invariants}. 
\end{proof}
\begin{remark} \label{rem:K'}
For parameter $\kappa \in (0,\,1)$, the Jacobian elliptic 
function $\mathrm{dn}(y\,|\,\kappa')$ 
has the real half-period 
\begin{equation} \label{eq:K'}
  \mathrm{K}' = \mathrm{F}(\tfrac{\pi}{2}\,|\,\kappa').
\end{equation}
\end{remark}
\begin{lemma} \label{th:cases}
Given a potential $\xi dz$ in normal form 
with $a,\,b \in \R$, then special choices 
for $a,\,b$ have the following effect: 
\begin{enumerate}
\item If $ab \neq 0$, then the surfaces generated by $\xi$ 
and $\xi^t$ differ by a rigid motion. Hence, swapping the roles 
of $a$ and $b$ gives the same surface up to isometry. 
\item The case $ab=0$ gives a twice punctured round sphere 
or a degenerate 'surface' consisting of a point. 
\item The case $a=\pm b \neq 0$ results in an 
associated family of a right circular cylinder. 
\end{enumerate}
\begin{proof} 
(i) Swapping the roles of the constants $a,\,b$ has no effect 
on the Hopf differential. We will show that the resulting 
surfaces have the same metric after a conformal change of 
coordinate (in fact a translation by a half-period). 
Let $a,\,b \in \R^*$, and let $v_1$ and $v_2$ 
be the conformal factors generated by  
$\xi dz$ respectively $\xi^t\,dz$, and $\kappa$ 
as in \eqref{eq:kappa}. 
Then $v_1(y) = 2b/H\,\mathrm{dn}(2by\,|\,1-\kappa)$, 
$v_2(y) = \sigma \, 2a/H\,\mathrm{dn}(2ay\,|\,1-\kappa^{-1})$ 
by Theorem \ref{th:invariants}, with the sign $\sigma = \pm 1$ 
in $v_2$ such that and $\sigma\,a$ and $b$ have the same sign. 
Then 
\begin{align*} 
   v_1(y + \tfrac{\mathrm{K}'}{2b}) &= 
   2bH^{-1}\,\mathrm{dn}(2by + \mathrm{K}'\,|\,\kappa') \\
   &= 2bH^{-1}\,\sqrt{1-\kappa'}\,\mathrm{nd}(2by\,|\,\kappa') 
   \mbox{ using 16.8.3 \cite{AbrS} } \\
   &= \sigma \, 2aH^{-1}\,\mathrm{dn}
     (2ay\,|\,-\tfrac{\kappa'}{1-\kappa'}) 
   \hfill \mbox{ using 16.10.4 \cite{AbrS} } \\
   &= \sigma \,2aH^{-1}\,\mathrm{dn}(2ay\,|\,1-\kappa^{-1}) 
     = v_2(y). 
\end{align*}
Thus after the change of coordinate 
$z \mapsto z + \tfrac{i\mathrm{K}'}{2|b|}$ the surfaces 
have same mean cuvature, Hopf differential and metric, and 
thus differ by a rigid motion. 
Hence swapping constants $a,\,b$ in $\xi dz$ 
gives the same surface up to a rigid motion.

\noindent (ii) In case $a=0$ and $b\neq 0$ the resulting 
immersion is either totally umbilic, or in the case of 
the parallel `surface' a point. In case $a=b=0$ 
the parallel surfaces both degenerate to a point, since 
then the correponding frame is $F \equiv \one$.

\noindent (iii) Since $\mathrm{dn}(u\,|\,0) \equiv 1$, 
the resulting conformal factor is $v(y) \equiv 2|b/H|$. 
Hence the associated family surface is that of a 
right circular cylinder.
\end{proof}
\end{lemma}
If we omitt the case $a = \pm b$, then as a 
consequence of Lemma \ref{th:cases} (i) we may assume 
without loss of generality that $a^2<b^2$ and thus that 
the parameter satisfies
\begin{equation}
  \kappa = \frac{a^2}{b^2} \in (0,\,1).
\end{equation}
We will need some properties of the square root $v_0$ of the 
Delaunay conformal factor.
\begin{proposition} \label{th:v_bounds}
For $\kappa \in (0,\,1)$, the function 
$v_0(y) = 2bH^{-1}\,\mathrm{dn}(2by\,|\,\kappa')$ 
satisfies
\begin{align}
  &v_0(y) = v_0(y + b^{-1}\mathrm{K}') 
  \mbox{ for all } y \in \R, \label{eq:v_period} \\
  &v_0(y + \tfrac{\mathrm{K}'}{2b}) =  
  v_0(\tfrac{\mathrm{K}'}{2b} - y) 
  \mbox{ for all } y \in \R, \label{eq:v_symmetry} \\
  &\dot{v}_0(y) = 0 \Longleftrightarrow  
  y \in \tfrac{\mathrm{K}'}{2b} \, \Z, \label{eq:v'is0} \\
  &v_0(y) = 2bH^{-1} \Longleftrightarrow  
  y \in b^{-1} \mathrm{K}'\, \Z. \label{eq:vmax}
\end{align}
\begin{proof}
Using Remark \ref{rem:K'} shows that $v_0$ has period 
$b^{-1}\mathrm{K}'$, proving \eqref{eq:v_period}. 
Assertion \eqref{eq:v_symmetry} 
follows from the change of argument 
formula 16.8.3 in \cite{AbrS}. 
Recall from Remark \ref{rem:v_elliptic} (i)  
that the derivative of $v_0$ vanishes at values 
$v_0 = \pm 2a/H,\, \pm 2b/H$, and with the above this proves 
\eqref{eq:v'is0} and \eqref{eq:vmax}.
\end{proof}
\end{proposition}


\section{The extended Delaunay frame}


To solve the period problem for a helicoidal {\sc{cmc}} surface,  
we need to compute the extended framing of an 
associated family of a Delaunay surface. This was obtained in 
\cite{KilKRS} in an untwisted setting, but since we need most 
of the ingredients, we provide the reader with all the details. 

Let $\xi dz$ be a potential in normal form with 
$a,\,b \in \R^*$ and $a \neq \pm b$, and $Q$ and 
$v_0$ as in Theorem \ref{th:invariants} with $\kappa < 1$. 
To find the unitary factor of the map $z \mapsto \exp(z\,\xi\,)$ 
it suffices to factorise $\exp(iy\,\xi\,)$, since 
$\exp(z\,\xi\,)$ is $\Lambda \SU_\tau$ 
valued along $y=0$. Off $y=0$, we make the following Ansatz: 
\begin{equation*} 
  F = \exp \left(\mathbf{f}\, \xi\, \right) T
\end{equation*}
where $\mathbf{f} = \mathbf{f}(x,\,y,\,\lambda) = 
x + iy + f(y,\,\lambda)$ for a function $f$ with 
$f(0,\,\lambda) = 0$, and $T = T(y,\,\lambda)$ 
an upper-triangular matrix 
\begin{equation*}
  T = \begin{pmatrix} A & B \\ 0 & A^{-1} \end{pmatrix},
\end{equation*}
with $A(0,\,\lambda) = 1$ and $B(0,\,\lambda) = 0$. 
The Maurer--Cartan form 
$\alpha = \alpha_1\,dx + \alpha_2\,dy$ of the 
unitary frame is given by
\begin{align} 
  \alpha_1 &= \frac{i}{2} \begin{pmatrix} \dot{v_0}\,v_0^{-1} & 
    2 v_0^{-1}\lambda^{-1}Q +  v_0\,\lambda \,H \\  
    v_0\lambda^{-1} H + 2v_0^{-1} \lambda\,\overline{Q} +  
    & - \dot{v_0}\,v_0^{-1} \end{pmatrix}, \label{eq:alpha_1} \\
  \alpha_2 &= \frac{1}{2} 
  \begin{pmatrix} 0 & - 2v_0^{-1}\lambda^{-1}Q + 
    v_0\,\lambda\, H \\ - v_0\,\lambda^{-1} H + 
    2v_0^{-1}\lambda\,\overline{Q} & 0 
  \end{pmatrix}. \label{eq:alpha_2}
\end{align}
For $F = \exp \left(\mathbf{f}\, \xi\, \right) T$ 
we then have
$F^{-1}dF = T^{-1}\,\xi\,\,T\,dx + 
  (\,\dot{\mathbf{f}} \,T^{-1}\,\xi\,\,T\ + 
  T^{-1}\,\dot{T}\,) \,dy$, 
so $\mathbf{f},\,T$ must simultaneously solve 
$\alpha_1 = T^{-1}\xi\,\,T$ and 
$\alpha_2 =  \dot{\mathbf{f}} \,
T^{-1}\xi\,\,T\ + T^{-1}\,\dot{T}$. 

We denote the lower left entry of $\alpha_1$ by $\Omega_1$, and 
since $Q = 2abH^{-1}$, we have
\begin{equation} \label{eq:Omeganz} \begin{split}
  \Omega_1(y,\,\lambda) &= \tfrac{i}{2} v_0(y)\,H\,\lambda^{-1} + 
  2 i\,v_0(y)^{-1}ab\,H^{-1}\,\lambda \\ 
  &\neq 0 \mbox{ for all } \lambda \in S^1,\,y \in \R. 
  \end{split}
\end{equation}
The claim $\Omega_1 \neq 0$ follows from 
remark \ref{rem:v_elliptic} (ii), 
and $a \neq \pm b$.

Let us denote the lower left entry of $\xi\,$ by $\omega$, and 
note that $a \neq \pm b$ ensures 
\begin{equation} \label{eq:omeganz}
   \omega =  ib\lambda^{-1} + ia\lambda \neq 0 \mbox{ for all } 
  \lambda \in S^1.
\end{equation}
Then $T^{-1} \xi\, \,T = \alpha_1$ is equivalent to 
\begin{equation*}
  \begin{pmatrix} 
    -A\,B\,\omega & -A^{-2}\omega^* - B^2 \omega \\ 
    A^2 \omega & A\,B\,\omega 
  \end{pmatrix} = 
  \begin{pmatrix} \frac{i}{2} v_0^{-1}\dot{v_0} & 
    -\Omega_1^* \\  
    \Omega_1  & -\frac{i}{2} v_0^{-1} \dot{v_0} 
  \end{pmatrix},
\end{equation*} 
with unique (up to sign) solution
\begin{equation} \label{eq:T}
  T = \frac{1}{\sqrt{\omega\,\Omega_1}}\,\begin{pmatrix}
  \Omega_1 & -\frac{i}{2}v_0^{-1}\dot{v_0} \\ 0 & \omega 
  \end{pmatrix}.
\end{equation}
Hence $T(y,\,\lambda)$ is defined for all $y \in \R$ and 
$\lambda \in S^1$. For $v_0$ as in \eqref{eq:v_Delaunay}, we have 
\begin{equation}
  \Omega_1(0,\,\lambda) = \omega(\lambda),
\end{equation}
and consequently $T(0,\,\lambda) = \one$. It remains to determine 
the function $\mathbf{f}$ such that 
\begin{equation} \label{eq:bf}
  \alpha_2 = \dot{\mathbf{f}} \,\alpha_1 + T^{-1}\dot{T}.
\end{equation}
To this end we compute
\begin{equation*}
  T^{-1}\dot{T} = \frac{1}{2\,\Omega_1}\begin{pmatrix} 
    \dot{\Omega}_1 & -i(\dot{v_0} v_0^{-1})\spdot \\ 
    0 & -\dot{\Omega}_1 \end{pmatrix} 
  = \frac{1}{2\,\Omega_1} \begin{pmatrix} 
    \dot{\Omega}_1 & iH^2v_0^2- 16iv_0^{-2}a^2b^2H^{-2} \\ 
    0 & -\dot{\Omega}_1 \end{pmatrix}.
\end{equation*}
We denote the lower left entry of $\alpha_2$ by 
\begin{equation} \label{eq:Omega2}
  \Omega_2 = -\tfrac{1}{2} v_0\,H\,\lambda^{-1} + 
  2 v_0^{-1}ab\,H^{-1}\,\lambda.
\end{equation}
The bottom left of \eqref{eq:bf} immediately yields 
$\dot{\mathbf{f}} = i + \dot{f} = \Omega_1^{-1}\Omega_2$. 
A quick computation shows that this solution  
is miraculously compatible with the other 
three equations in \eqref{eq:bf}. Thus 
$\dot{f} = \Omega_1^{-1}(\Omega_2 - i \Omega_1)$, and 
we define
\begin{equation} \label{eq:f} 
  f(y,\,\lambda) = \int_0^y \frac{\Omega_2(\zeta,\,\lambda) - 
    i \Omega_1(\zeta,\,\lambda)}
    {\Omega_1(\zeta,\,\lambda)}\,d\zeta 
  = \int_0^y \frac{8ab\lambda\,d\zeta}{4iab\lambda + 
    iH^2\lambda^{-1}v_0^2(\zeta)},
\end{equation}
and now set 
\begin{equation} \label{eq:bold_f}
  \mathbf{f}(x,\,y,\,\lambda) = x + iy + f(y,\,\lambda).
\end{equation}
Clearly $f(0,\,\lambda) = 0$ and we thus have a solution of 
$dF = F\,\alpha_1\,dx + F\,\alpha_2\,dy$ 
with $\alpha_1,\,\alpha_2$ as in \eqref{eq:alpha_1} 
and \eqref{eq:alpha_2} of the form 
$F=\exp(\mathbf{f}\,\xi\,)\,T$. 

For $v_0$ as in \eqref{eq:v_Delaunay}, 
it is easy to verify that $T(0,\,\lambda) = \one$, and 
consequently $F(0,\,0,\,\lambda) = \one$. 
To show that $F$ is the extended framing of $\xi\,$, 
it remains to show that 
$B(y,\,\lambda) = T^{-1} \exp(\,-f\,\xi\, \,) 
\in \Lambda_+\SL_\tau$ for all $y \in \R$. 
From $\Phi = F B$, it follows that $B = B(y,\,\lambda)$ 
solves $dB + \alpha\,B = B\,\xi\,\,dz$ with 
$B(0,\,\lambda) = \mathrm{Id}$. This decouples into 
$\dot{B} + \alpha_2\,B = i\,B\,\xi\,$ and 
$\alpha_1\,B = B\,\xi\,$. 
It thus remains to verify that
\begin{equation} \label{eq:B_ode} 
 \dot{B}\,B^{-1} = i\,\alpha_1 - \alpha_2.
\end{equation}
On the one hand, $\dot{B}\,B^{-1} = -T^{-1}\dot{T} - 
\dot{f}\,\alpha_1$, while \eqref{eq:alpha_1} and 
\eqref{eq:alpha_2} yield
\begin{equation*}
  i\,\alpha_1 - \alpha_2 = \frac{1}{2} 
    \begin{pmatrix} -\dot{v_0}\,v_0^{-1} & 
    -2v_0\,\lambda \,H \\ -4v_0^{-1}\lambda\,\bar{Q} & \dot{v_0}\,v_0^{-1} 
  \end{pmatrix}.
\end{equation*}
A computation verifies that indeed $B$ 
solves \eqref{eq:B_ode} and hence is positive 
and smooth. Further, \eqref{eq:B_ode} can be 
integrated at $\lambda = 0$ and gives 
$B_0 = B(y,\,0) = \mathrm{diag}[\,v_0^{-1/2},\,v_0^{1/2}\,]$. 
Hence $B:\R \to \Lambda_+\SL_\tau$. Altogether, this proves
\begin{theorem}\cite{KilKRS} 
The extended framing for the Delaunay surface 
generated by a potential $\xi dz$ in normal form 
with $a \neq \pm b$ is given by 
\begin{equation}\label{eq:Delaunay_frame}
  F = \exp \left(\mathbf{f}\, \xi\, \right) T
\end{equation}
with $\mathbf{f} = \mathbf{f}(x,\,y,\,\lambda)$ 
as in \eqref{eq:bold_f} and $T = T(y,\,\lambda)$ 
as in \eqref{eq:T} and $v_0$ as in \eqref{eq:v_Delaunay}.
\end{theorem} 


\section{Helicoidal {\sc{cmc}} cylinders in $\R^3$}


%
We now arrive at our main objective which is to solve period 
problems for equivariant harmonic maps $\C \to S^2$ as well 
as for the corresponding constant mean curvature surfaces. 
We will show that every equivariant harmonic map $\C \to S^2$ 
is periodic. Further, we will discuss the closing conditions 
for helicoidal constant mean curvature cylinders. 

Other methods have been used to study helicoidal {\sc{cmc}} 
surfaces, as in the investigations of Roussos et al., 
see \cite{HitR} and \cite{Rou_heli} and the references therein.  

Consider translations 
\begin{equation*}
  \delta: z \mapsto z + p + iq \mbox{ with } p,\,q \in \R\, .
\end{equation*}
If an extended framing $F(z,\,\lambda)$ with 
$F(0,\,\lambda) = \one$ has a well defined 
monodromy 
\begin{equation}
  M(\delta,\,\lambda) = \delta^*F\,F^{-1} 
  = F(\delta(0),\,\lambda)
\end{equation}
with respect to $\delta$, 
then $\delta^*f_{\mu_0} = f_{\mu_0}$ for $\mu_0 \in S^1$ 
if and only if $\Phi_{\mu_0}( M ) = \one$, 
which in turn is equivalent to 
the two \emph{closing conditions}:
\begin{equation} \label{eq:closing}
  M(\delta,\,\mu_0) = \pm \one \mbox{ and } 
  M'(\delta,\,\mu_0) = 0.
\end{equation}
We first consider the case of a right circular cylinder, 
generated by a potential in normal form with $a=\pm b \neq 0$. 
\begin{lemma}
Every member in the associated family of a right circular 
cylinder has a period.
\end{lemma}
\begin{proof}
Up to rigid motion, an associated family $f_\lambda$ of a 
right circular cylinder is generated by 
a potential $\xi dz$ in normal form with $a= \pm b$. 
Consider the case $a=b$. 
The corresponding extended frame is 
$\exp(ia(z\lambda^{-1} + \bar{z}\lambda) 
  \bigl(\begin{smallmatrix}0&1\\1&0\end{smallmatrix}\bigr) )$
with monodromy 
\begin{equation*}
  M(\delta,\,\lambda) = \cosh(ia(\delta\,\lambda^{-1} + 
  \bar{\delta}\,\lambda))\,\one
  + \sinh(ia(\delta\,\lambda^{-1} + 
  \bar{\delta}\,\lambda))\,
  \bigl( \begin{smallmatrix} 0 & 1 \\ 
    1 & 0 \end{smallmatrix} \bigr).
\end{equation*}
Hence $M(\delta,\,\lambda_0)$ satisfies the two 
closing conditions 
if and only if $\delta\lambda_0^{-1} + \bar{\delta}\lambda_0 \in 
\tfrac{\pi}{a}\Z$ and 
$\bar{\delta}\lambda_0 - \delta \lambda_0^{-1} = 0$, 
and thus if and only if $\delta\in\tfrac{\pi}{2a}\lambda_0\Z$. 
Then for all $k \in \Z$
\[
  f_{\lambda_0}(z+\tfrac{\pi}{2a}\lambda_0\,k) 
  = f_{\lambda_0}(z).
\] 
The case $a=-b$ is proven similarly. 
\end{proof}
Let $\xi\,$ as in \eqref{eq:del_mod} with $a,\,b \in \R^*$. 
It is easily verified that 
$0 \leq \det \xi\,(\lambda)$ for $\lambda \in S^1$, since 
$\min\{(a-b)^2,\,(a+b)^2\}$ is a lower bound. Hence 
\begin{equation} \label{eq:Delta}
  \mu(\lambda) = \sqrt{-\det \xi\,(\lambda)}.
\end{equation}
is purely imaginary for all $\lambda \in S^1$. 
The roots of $\mu$ are 
$\pm i \,\sqrt{a/b},\,\pm i \,\sqrt{b/a}$. 
Since we have discussed the case $a = \pm b$ of a 
right circular cylinder, 
we will from now assume that $a \neq \pm b$ and thus ensure 
that these roots are off $S^1$. 

To study the monodromy behaviour of extended frames 
generated by potentials in normal form, define 
\begin{equation} \label{eq:bold_g}
  \mathbf{g} = \mathbf{g}(\delta,\,\lambda) = 
  \mathbf{f}(p,\,q,\,\lambda)\,\mu(\lambda),
\end{equation}
with $\mathbf{f}$ as in \eqref{eq:bold_f} and 
$\mu(\lambda)$ as in \eqref{eq:Delta}. 
The monodromy of $F$ in \eqref{eq:Delaunay_frame}
with respect to $\delta$ is then 
$M(\delta,\,\lambda) = F(p,\,q,\,\lambda) = 
  \exp\left( \mathbf{f}(p,\,q,\,\lambda)\,
    \xi\,(\lambda) \right)\,
  T(q,\,\lambda)$, which we may write as 
\begin{equation} \label{eq:Delaunay_monodromy} 
  M(\delta,\,\lambda) = \cosh \left(\mathbf{g}\right)\,
  T  + 
  \mu^{-1}
  \sinh \left(\mathbf{g}\right)\,
  \xi\,T. 
\end{equation}
We have seen in Lemma \ref{th:normalform} that at 
$\lambda^4 = 1$ we obtain surfaces of revolution. 
We omit these fourth roots of unity in 
our further study of helicoidal cylinders. 
\begin{lemma} \label{th:closing_conditions}
Let $M(\delta,\,\lambda)$ as in \eqref{eq:Delaunay_monodromy} 
be the monodromy of an extended frame generated by a potential 
in normal form with $a \neq \pm b$, and 
$\lambda_0 \in S^1\setminus\{\pm 1,\,\pm i \}$.

\noindent {\rm{(i)}} $M(\delta,\,\lambda_0) = \pm \one$ 
if and only if $T(q,\,\lambda_0) = \one$ and 
\begin{equation}
  \mathbf{g}(\delta,\,\lambda_0) \in \pi i \Z, \label{eq:g} 
\end{equation}

\noindent {\rm{(ii)}} $T(q,\,\lambda_0) = \one$ 
holds if and only if $v_0(q) = v_0(0)$ and $\dot{v_0}(q) = 0$. 

\noindent {\rm{(iii)}} If we assume 
$M(\delta,\,\lambda_0) = \pm \one$, 
then $M'(\delta,\,\lambda_0) = 0$ if and only if 
  \begin{equation}
    \mathbf{g}'(\delta,\,\lambda_0) = 0. \label{eq:g'}
  \end{equation}
\end{lemma}
\begin{proof}
(i) If $M(\delta,\,\lambda_0) = \pm \one$, then the lower left 
entry in \eqref{eq:Delaunay_monodromy} is 
\begin{equation*}
  \mu(\lambda_0)^{-1} \sinh(\mathbf{g}(\lambda_0)) \,
  \sqrt{\omega(\lambda_0) \,\Omega_1(q,\,\lambda_0)} = 0.
\end{equation*}
Hence either $\mathbf{g}(\delta,\,\lambda_0) \in \pi i \Z$ or 
$\omega(\lambda_0) \,\Omega_1(q,\,\lambda_0) = 0$. 
But the latter is not possible by \eqref{eq:Omeganz} and 
\eqref{eq:omeganz}. 
Hence $\mathbf{g}(\delta,\,\lambda_0) \in \pi i \Z$, 
and consequently $T(q,\,\lambda_0) = \one$. 
The converse is obvious. 

\noindent (ii) Recall the definition of $T$ from 
equation \eqref{eq:T}. 
Clearly $T(q,\,\lambda_0) = \one$ holds if and only if 
$\omega(\lambda_0) = \Omega_1(q,\,\lambda_0)$ and 
$\dot{v_0}(q) = 0$. 
Now $\omega(\lambda_0) = \Omega_1(q,\,\lambda_0)$ is 
a quadratic equation in $v_0(q)$ with solutions 
$v_0(q) = 2b/H,\,\lambda_0^2/(5H)$. Since $\dot{v_0}(q) = 0$, 
we must have $v_0(q) = 2b/H = v_0(0)$ by \eqref{eq:v_elliptic}.

\noindent (iii) Differentiating the monodromy in 
\eqref{eq:Delaunay_monodromy} with respect to $t$, 
and making use of  
$T(q,\,\lambda_0) = \one$ and 
$\mathbf{g}(\delta,\,\lambda_0) \in \pi i \Z$, 
we obtain
\begin{equation} \label{eq:M'}
  M'(\delta,\,\lambda_0) = \pm T'(q,\,\lambda_0) \pm 
  \mu(\lambda_0)^{-1}
  \mathbf{g}'(\delta,\,\lambda_0)\,\xi\,(\lambda_0) .
\end{equation}
In part (ii) we saw that $T(q,\,\lambda_0) = \one$ 
implies that $\dot{v_0}(q) = 0$ and 
$\omega(\lambda_0) = \Omega_1(q,\,\lambda_0)$. 
Thus the first term on the right hand side in \eqref{eq:M'} is 
diagonal and 
\begin{equation*}
  T'(q,\,\lambda_0) = 
    \tfrac{1}{2\omega}(\Omega_1'(q,\,\lambda_0) - 
    \omega'(\lambda_0))\,\,\one. 
\end{equation*}
On the other hand, the second term on the right hand side 
in \eqref{eq:M'} is off-diagonal, and since the entries of 
$\xi\,$ have no roots on $S^1$, the vanishing of 
$M'(\delta,\,\lambda_0)$ is 
equivalent to the simultaneous vanishing of 
$\mathbf{g}'(\delta,\,\lambda_0)$ and 
$(\Omega_1'(q,\,\lambda_0) - \omega'(\lambda_0))$.  
A quick computation reveals that $v_0(q) = 2b/H$ ensures 
that $\Omega_1'(q,\,\lambda_0) = \omega'(\lambda_0)$ holds. 
Hence $M'(\delta,\,\lambda_0) = 
\pm \,\mathbf{g}'(\delta,\,\lambda_0)\,
\xi\,(\lambda_0)/\mu(\lambda_0)$ and 
the claim follows. 
\end{proof}

 
\subsection{The rotational period} 
Solving a period problem for a constant mean curvature surface 
is tantamount to ensuring the two closing 
conditions \eqref{eq:closing}. The closing conditions are 
equivalent to killing rotational 
and translational periods respectively. We first consider 
rotational periods, and will learn that one of the generators 
of the period lattice has to be a period of the metric. 
We also show that every equivariant harmonic map $\C \to S^2$ 
is periodic.

Consider the function $v_0$ from \eqref{eq:v_Delaunay}. 
We have seen in \eqref{eq:v'is0} and \eqref{eq:vmax} that 
the conditions $\dot{v_0}(q) = 0$ and $v_0(q) = v_0(0)$ hold if 
and only if $q \in b^{-1} \mathrm{K}'\,\Z$. So let 
\begin{equation} \label{eq:q}
  q =  |b|^{-1} \mathrm{K}' k \,\mbox{ for some } \,k \in \Z.
\end{equation}
For any such $q$ we have $v_0(y+q) = v_0(y)$ and 
$T(y+q,\,\lambda) = T(y,\,\lambda)$. Consequently,
\begin{equation*} \begin{split}
  f(y+q) &= \int_0^{y+q} \frac{8ab\lambda\,d\zeta}{4iab\lambda + 
    iH^2\lambda^{-1}v_0^2(\zeta)} \\
  &= \int_0^y \frac{8ab\lambda\,d\zeta}{4iab\lambda + 
    iH^2\lambda^{-1}v_0^2(\zeta)} + 
  \int_0^q \frac{8ab\lambda\,d\zeta}{4iab\lambda + 
    iH^2\lambda^{-1}v_0^2(\zeta)} \\ 
  &= f(y) + f(q), \end{split}
\end{equation*}
and therefore
\begin{equation}
  F(x+p,\,y+q,\,\lambda) = 
  \exp((p+iq+f(q))\xi)\,F(x,\,y,\,\lambda).
\end{equation}
Since $F$ takes values in $\Lambda \SU_\tau$, it follows that 
\begin{equation} \label{eq:Imf}
  \mathrm{Im}(f(q,\,\lambda)) = 
  -q \mbox{ for all } \lambda \in S^1.
\end{equation}
We now proceed to find values for 
$p \in \R$ such that for some fixed 
$\lambda_0 \in S^1 \setminus \{ \pm 1,\,\pm i \}$ 
both \eqref{eq:g} and \eqref{eq:g'} hold. 

Recall the definitions of $\mathbf{g}(\delta,\,\lambda_0)$ 
from \eqref{eq:bold_g} 
and $f(q,\,\lambda_0)$ from \eqref{eq:f}. 
Now $\mathbf{g}(\delta,\,\lambda_0) \in \pi i \Z$ 
in \eqref{eq:g} means that 
$\left( p + i\,q + f(q,\,\lambda_0)\right)
  \mu(\lambda_0) = \pi \,i\,l$ for some $l \in \Z$. 
Solving this for $p$, and using \eqref{eq:Imf}, we have that 
$\mathbf{g}(\delta,\,\lambda_0) \in \pi i \Z$ if and only if 
$p$ is of the form $p  = \pi\,i\, l \,\mu^{-1}(\lambda_0) - 
  \mathrm{Re}(f(q,\,\lambda_0))$ for some $l \in \Z$. 
We have proven
\begin{lemma}
The monodromy $M(\delta,\,\lambda)$ of an extended frame 
generated by a potential in normal form satisfies 
$M(\delta,\,\lambda_0) = \pm \one$ if and only if 
$\delta: z \mapsto z+p+iq$ with 
\begin{align}
  p  = \pi\,i\, l \,\mu^{-1}(\lambda_0) - 
  \mathrm{Re}(f(q,\,\lambda_0)) \mbox{ for some } l \in \Z&, 
  \label{eq:p_simple} \\
  q = |b|^{-1}\mathrm{K} \,k \mbox{ for some } k \in \Z&. 
  \label{eq:q_simple}
\end{align}
\end{lemma}
\begin{theorem}
Every equivariant harmonic map $N:\C \to S^2$ is periodic.
\end{theorem}
\begin{proof}
By Corollary \ref{cor:equi_harm}, up to isometry, 
every associated family of equivariant harmonic maps 
$N_\lambda:\C \to S^2$ can be generated by a potential $\xi dz$
in normal form. Assume, possibly after a rigid motion, that 
$N = N_{\lambda_0}$ for some $\lambda_0 \in S^1$. 
The extended framing then has monodromy $M(\delta,\,\lambda)$ 
as in \eqref{eq:Delaunay_monodromy}. The periodicity 
of the harmonic map $N(z+\delta) = N(z)$ 
is equivalent to the first closing 
condition $M(\delta,\,\lambda_0) = \one$, which holds 
for $\delta = p + iq$ with $p$ as in \eqref{eq:p_simple} and 
$q$ as in \eqref{eq:q_simple}.
\end{proof}
%


\subsection{The translational period} 
Amongst equivariant harmonic maps, we now single out 
those which generate periodic immersions. 

It remains to realize the constraint of \eqref{eq:g'}, 
as a condition on the spectral parameter. For 
$\lambda_0 = e^{it_0} \in S^1$ and 
$p,\,q$ as in \eqref{eq:p_simple} and \eqref{eq:q_simple} 
we have that $(\mathbf{f}\,\mu)(p,\,q,\,\lambda_0) = \pi i l$ 
for some $l \in \Z$, and compute 
\begin{align*}
  \mathbf{g}'(p,\,q,\,\lambda_0) &= 
  f'(q,\,\lambda_0)\,\mu(\lambda_0) +  
  \mathbf{f}(p,\,q,\,\lambda_0)
  \,\mu'(\lambda_0)\\ 
  &=  f'(q,\,\lambda_0)\,\mu(\lambda_0) + 
  \pi i l\,\mu^{-1}(\lambda_0) 
  \,\mu'(\lambda_0).
\end{align*}
Hence $\mathbf{g}'(p,\,q,\,\lambda_0) = 0$ is 
equivalent to 
\begin{equation} \label{eq:l}
  l = - \frac{\mu^2(\lambda_0)}
  {\pi\,i \,\mu'(\lambda_0)}\,f'(q,\,\lambda_0).
\end{equation}
Since it is not apparent for which $\lambda \in S^1$, if 
indeed any, the right hand side in \eqref{eq:l} 
is integer-valued, 
we view it as a function of $\lambda$ and define 
\begin{equation} \label{eq:L}
  L(\lambda) = 
  \,- \frac{\mu^2(\lambda)}
  {\pi\,i \,\mu'(\lambda)}
  f'(q,\,\lambda).
\end{equation}
\begin{proposition}
The function $L(\lambda)$ is real and meromorphic on 
$S^1$ with only simple poles at the fourth roots of unity.
\begin{proof} 
Consider first the function $f':S^1 \to \R^*$ given by 
\begin{equation} \label{eq:f'}
  f'(q,\,\lambda) = \int_0^q \frac{16ab\,H^2\,
    v_0^2(\zeta)\,\,d\zeta}
    {(4ab\lambda + H^2\lambda^{-1}v_0^2(\zeta))^2}.
\end{equation}
By \eqref{eq:Omeganz}, 
the denominator in the integrand 
in \eqref{eq:f'} is non-zero for all $\lambda \in S^1$. 
By \eqref{eq:Imf} we have that 
$f'(q,\,\lambda) = \mathrm{Re}f'(q,\,\lambda)$ 
for each $\lambda \in S^1$, 
and conclude that $f'$ is real analytic and non-zero. 
The first factor of $L(\lambda)$ computes to
\begin{equation} \label{eq:2} 
  - \frac{\mu^2(\lambda)}{\pi\,i \,\mu'(\lambda)} = 
  \frac{\mu^3(\lambda)}{\pi ab\,(\lambda^{-2}-
    \lambda^2)},
\end{equation}
which is clearly real for $\lambda \in S^1$, since 
$\mu(\lambda) \in i \R^*$ and $\lambda^{-2}-\lambda^2 
\in i\R$. Hence $L(\lambda)$ is real for $\lambda \in S^1$. 
The only poles of $L$ are those of the first factor 
\eqref{eq:2}, namely simple poles at the fourth roots 
of unity. This conludes the proof.
\end{proof}
\end{proposition}

As a consequence of the last proposition, we are assured 
of infinitely many values 
$\lambda \in S^1\setminus \{ \lambda^4=1\}$ 
at which $L(\lambda) \in \Z$. 
This proves the existence of infinitely many helicoidal 
constant mean curvature cylinders, and we conclude  
\begin{theorem}
In each associated family of a Delaunay surface there are 
infinitely many non-congruent cylinders with screw-motion 
symmetry.
\end{theorem}


\bibliographystyle{amsplain}

\providecommand{\bysame}{\leavevmode\hbox to3em{\hrulefill}\thinspace}
\providecommand{\MR}{\relax\ifhmode\unskip\space\fi MR }
\providecommand{\MRhref}[2]{%
  \href{http://www.ams.org/mathscinet-getitem?mr=#1}{#2}
}
\providecommand{\href}[2]{#2}

\end{document}